# The Generalized *RS*-theorem:

## Bounds for the second largest eigenvalue of a connected graph with a cut-vertex

Bojana Mihailović, Marija Rašajski

*Faculty of Electrical Engineering, University of Belgrade*

## Abstract

Eigenvalues of a graph are the eigenvalues of the corresponding *(0,1)*-adjacency matrix. The second largest eigenvalue ($\lambda_2$) provides significant information on characteristics and structure of graphs. Therefore, finding bounds for $\lambda_2$ is a topic of interest in many fields. So far we have studied the graphs with the property $\lambda_2 \leq 2$, so-called *reflexive* graphs. The original *RS*-theorem is about them. In this paper we generalize that concept and introduce the arbitrary bounds. The Generalized *RS*-theorem gives us an answer whether the second largest eigenvalue of a graph is greater than, less than, or equal to *a*, *a>0*, within some classes of connected graphs with a cut-vertex. After removing the cut-vertex *u* of the given graph *G*, we examine the indices of the components of *G-u*. The information on these indices is used to make conclusions about the second largest eigenvalue of the graph *G*. In the *RS*-theorem *a=2*. Here, we state and prove the Generalized *RS*-theorem, one corollary, and, also, give some useful lemmas.

## 1. Introduction

We consider simple connected graphs, undirected, with no loops or multiple edges. If $A$ is the $(0,1)$-adjacency matrix of the graph $G$, then its *characteristic polynomial* is defined by $P_G(\lambda) = \det(\lambda I - A)$. The roots of the characteristic polynomial are called *eigenvalues* of $G$. The family of eigenvalues is the *spectrum* of $G$. $A$ is a real and symmetric matrix, and, therefore, its eigenvalues are real. We assume their non-increasing order: $\lambda_1(G) \geq \lambda_2(G) \geq ... \geq \lambda_n(G)$. The largest eigenvalue $\lambda_1(G)$ is called the *index* of $G$. For connected graphs $\lambda_1(G) > \lambda_2(G)$ holds. If graph $G$ is not connected, then its spectrum consists of all the eigenvalues of its components (with all their multiplicities). If $G$ is disconnected, and if the two largest indices of its components are equal, then $\lambda_1(G) = \lambda_2(G)$. In this paper the when we use the term *subgraph*, it means the *induced subgraph*.

The second largest eigenvalue of a graph is a subject of investigations in spectral graph theory, but also in computer science and various fields across the science in which networks as mathematical models are widely used. For instance, graphs for which $\lambda_2 \leq 2$, also called *reflexive graphs*, correspond to sets of vectors in the Lorentz space $R^{p,1}$, and they are Lorentzian counterparts of the spherical and Euclidean graphs. They have direct application in the construction and the classification of reflection groups [10].

Also, the second largest eigenvalue is important in determining the structure of connected *regular* graphs, graphs whose all vertices have the same degree. For example, for cubic graphs, regular graphs of degree *3*, it was observed that for small values of $\lambda_2$ graphs have smaller diameters, higher connectivity and girth, and for large values of $\lambda_2$, they have larger diameter, lower connectivity and girth [4].

In spectral graph theory, graphs with the second largest eigenvalue bounded by a constant $a \in R$ have been investigated by many authors. Some of the bounds considered so far are: $a = 1/3$ [1], $a = \sqrt{2} - 1$ [11], $a = (\sqrt{5} - 1)/2$ [5, 21], $a = 1$ [2, 9]. *Reflexive graphs*, have been investigated in [6, 7, 8, 9, 12, 13, 14, 15, 16, 17, 18, 19].

The *RS*-theorem [17] can often be used to prove whether the connected graph *G* with a cut-vertex *u* is reflexive or not. After the removal of the cut-vertex we compare the components of *G-u* to Smith graphs (graphs with the index *2*), and draw conclusions from the *RS*-theorem.

The Generalized *RS*-theorem tells us whether $\lambda_2(G) < a$, $\lambda_2(G) > a$ or $\lambda_2(G) = a$, for some $a > 0$, within some classes of connected graphs with a cut-vertex *u*. The information on indices of the components of the graph *G-u*, namely how they are compared to *a*, is used to make conclusions about the second largest eigenvalue of *G*. In the *RS*-theorem $a = 2$.

This paper is structured in the following way. After the Introduction, we will present the main tools used in our investigations: the Interlacing Theorem and the Schwenk's Lemma. Then, we will show the Smith graphs and the *RS*-theorem. The main section of the paper brings the Generalized *RS*-theorem, along with four auxiliary lemmas that are used to prove the main theorem and a corollary. We predict that these lemmas will be useful practical tools in our future research. After that we show the cases when the Generalized *RS*-theorem in not applicable, and, also, we give so examples of appliacations.

## 2. Auxiliary results and the *RS*-theorem

In this section we will show the main tools that are used in our investigations, namely the Interlacing Theorem and the Schwenk's Lemma. Then, we will present the Smith graphs and the *RS*-theorem.

The following theorem shows the interrelation between the spectra of a graph and its induced subgraph.

**The Interlacing Theorem.** *Let* $\lambda_1 \geq \lambda_2 \geq ... \geq \lambda_n$ *be the eigenvalues of a simple graph G and* $\mu_1 \geq \mu_2 \geq ... \geq \mu_m$ *the eigenvalues of its induced subgraph H. Then the inequalities* $\lambda_{n-m+i} \leq \mu_i \leq \lambda_i \ (i=1,2,...,m)$ *hold.*

If $G$ is a connected graph and $m = n-1$, then $\lambda_1 > \mu_1 \geq \lambda_2 \geq \mu_2 \geq ...$.

**Schwenk's Lemma** (A. J. Schwenk, [20]) Given a graph $G$, let $C(v)$ and $C(uv)$ denote the set of all cycles containing a vertex $v$ and an edge $uv$ of $G$, respectively. Then

1) $P_G(\lambda) = \lambda P_{G-v}(\lambda) - \sum_{u \in Adj(v)} P_{G-v-u}(\lambda) - 2 \sum_{C \in C(v)} P_{G-V(C)}(\lambda),$

2) $P_G(\lambda) = P_{G-uv}(\lambda) - P_{G-v-u}(\lambda) - 2 \sum_{C \in C(uv)} P_{G-V(C)}(\lambda),$

where *Adj(v)* denotes the set of neighbours of *v*, while *G - V(C)* is the graph obtained from *G* by removing the vertices belonging to the cycle *C*.

The only connected graphs for which $\lambda_1 = 2$ holds are called *Smith graphs* [22] and they are shown in Figure 1. Proper subgraphs of Smith graphs all have the property $\lambda_1 < 2$ and they are called Coxeter-Dynkin graphs.

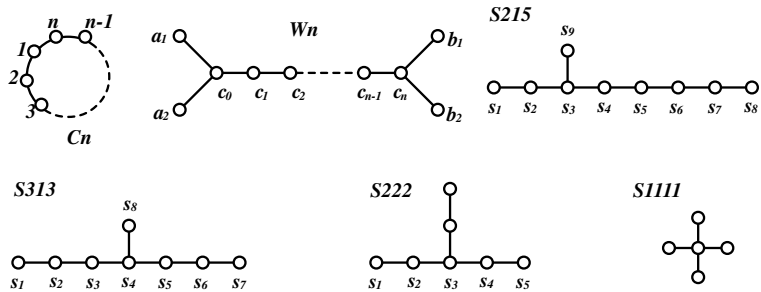

**Figure 1**

Whether the given graph with a cut-vertex is reflexive or not in many cases can be established effortlessly using the *RS*-theorem.

**Theorem RS.** [17] Let $G$ be a graph with a cut-vertex $u$.

1) If at least two components of $G$-$u$ are supergraphs of Smith graphs, and if at least one of them is a proper supergraph, then $\lambda_2(G) > 2$.

2) If at least two components of $G$-$u$ are Smith graphs, and the rest are subgraphs of Smith graphs, then $\lambda_2(G) = 2$.

3) If at most one component of $G$-$u$ is a Smith graph, and the rest are proper subgraphs of Smith graphs, then $\lambda_2(G) < 2$.

After removing a cut-vertex of the graph G, we get several new connected graphs. They are comparable to the Smith graphs, in the sense that they are either their subgraphs or supergraphs, and then we apply the *RS*-theorem.

In many cases the theorem gives the answer about the reflexivity of the graph, but there is one case when it does not. This is the case when after the removal of the cut-vertex $u$ we get one proper supergraph, and the rest are proper subgraphs of Smith graphs. In these cases other techniques are used, but often, at the subgraph level analysis, the *RS*-theorem proves useful again.

In previous work, we have described many classes of reflexive graphs [6-9, 13-19]. Since reflexivity is a hereditary property, every subgraph preserves it, then it is natural to present classes of reflexive graphs through the sets of maximal reflexive graphs. In this case maximal means that its supergraphs are not reflexive. Our research shows that Smith graphs play an essential role in the construction of the maximal reflexive graphs [13, 18, 19].

This is how we got the motivation to generalize these ideas, and among them the *RS*-theorem. We have a set of graphs with the property $\lambda_1 = 2$, and we have found various ways to construct the reflexive graphs (graphs with the property $\lambda_2 \leq 2$) using them. Now, we try to do something similar with the arbitrary bounds. Say, we know a set of graphs for which $\lambda_1 = a$, can we construct some classes of graphs with the property $\lambda_2 \leq a$? The first step towards answering this, and many other similar questions, is the generalization of the *RS*-theorem.

## 3. The Generalized RS-theorem

In this section we will state and prove the Generalized *RS*-theorem. Before that, we present four lemmas that will be used in the proof of the theorem.

**Lemma 1.** Consider graph $G$ in Figure 2, where $G_1$ is a connected graph with the index $a$, $a > 0$, and $u$ is an extra vertex connected to some of the vertices of the graph $G_1$ ($v_1, v_2, ..., v_m$). Then, $P_G(a) < 0$, and, consequently, $\lambda_2(G) < a < \lambda_1(G)$.

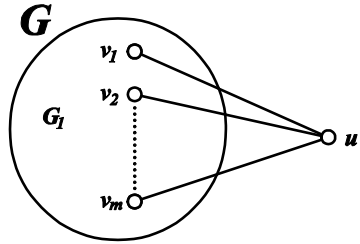

**Figure 2**

**Proof.** By the Interlacing theorem $\lambda_1(G) > a$ and $\lambda_2(G) \leq a$. Applying the Schwenk lemma at vertex $u$, we get the characteristic polynomial of the graph $G$:

$$P_G(\lambda) = \lambda P_{G_1}(\lambda) - \sum_{i=1}^{m} P_{G_1 - v_i}(\lambda) - 2 \sum_{C \in C(u)} P_{G-C}(\lambda).$$

$\lambda_1(G_1 - v_i) < a$ holds, implying $P_{G_1 - v_i}(a) > 0$ for $i = 1, ..., m$. $P_{G-C}(a) > 0$ holds, since graph $G - C$ ($C \in C(u)$) is a subgraph of the graph $G_1 - v_i$ for some $i = 1, ..., m$. Therefore, $P_{G_1}(a) = 0$ implies $P_G(a) < 0$ and $\lambda_2(G) < a$. □

The following two lemmas tell us more about the multiplicities of the eigenvalues in case we add or delete a vertex of an existing connected graph. They are the direct consequences of the Interlacing theorem.

**Lemma 2.** Let $G_1$ be a connected graph shown in Figure 2, and let $\lambda_1(G_1) > \lambda_2(G_1) \geq \lambda_3(G_1) \geq ... \geq \lambda_{n-1}(G_1) \geq \lambda_n(G_1)$ be its eigenvalues. Let $\alpha$ be an eigenvalue of the multiplicity $k$ such that $\lambda_{p+1}(G_1) = \lambda_{p+2}(G_1) = ... = \lambda_{p+k}(G_1) = \alpha$. Let $G$ be a graph obtained by adding a vertex to the graph $G_1$ and connecting it to some vertices of $G_1$. Then, $\lambda_{p+2}(G) = \lambda_{p+3}(G) = ... = \lambda_{p+k}(G) = \alpha$. Additionally, it is possible that $\lambda_{p+1}(G) = \alpha$ or $\lambda_{p+k+1}(G) = \alpha$. Therefore, the multiplicity of the eigenvalue $\alpha$ of graph $G$ is $k$, $k-1$, or $k+1$.

**Lemma 3.** Let $G$ be a connected graph shown in Figure 2, and let $\lambda_1(G) > \lambda_2(G) \geq \lambda_3(G) \geq ... \geq \lambda_{n-1}(G) \geq \lambda_n(G)$ be its eigenvalues. Let $\alpha$ be an eigenvalue of the multiplicity $k$ such that $\lambda_{p+1}(G) = \lambda_{p+2}(G) = ... = \lambda_{p+k}(G) = \alpha$. Let $G_1$ be a connected graph obtained by deleting a vertex of the graph $G$. Then, $\lambda_{p+1}(G_1) = \lambda_{p+2}(G_1) = ... = \lambda_{p+k-1}(G_1) = \alpha$. Additionally, it is possible that $\lambda_p(G_1) = \alpha$ (if $p > 1$) or $\lambda_{p+k}(G_1) = \alpha$ (if $p + k < n$). Therefore, the multiplicity of the eigenvalue $\alpha$ of graph $G_1$ is $k$, $k-1$, or $k+1$.

**Lemma 4.** Let $G$ be a graph with $n$ vertices and the eigenvalues $\lambda_1(G) \geq \lambda_2(G) \geq \lambda_3(G) \geq ... \geq \lambda_{n-1}(G) \geq \lambda_n(G)$. Let $\alpha = \lambda_{m+1}(G) = \lambda_{m+2}(G) = ... = \lambda_{m+k}(G)$ be the eigenvalue of multiplicity $k$. If the polynomial $Q_G(\lambda)$ is defined by the relation $P_G(\lambda) = (\lambda - \alpha)^k Q_G(\lambda)$, then $\text{sgn}(Q_G(\alpha)) = (-1)^m$.

**Proof.** The characteristic polynomial of the graph $G$ is factored in the following way:
$$P_G(\lambda) = \prod_{i=1}^{m}(\lambda - \lambda_i) \cdot (\lambda - \alpha)^k \cdot \prod_{i=m+k+1}^{n}(\lambda - \lambda_i).$$ Then, $Q_G(\lambda) = \prod_{i=1}^{m}(\lambda - \lambda_i) \cdot \prod_{i=m+k+1}^{n}(\lambda - \lambda_i)$. For $i \in \{1, 2, ..., m\}$ $\alpha - \lambda_i < 0$, and $\text{sgn}\left(\prod_{i=1}^{m}(\alpha - \lambda_i)\right) = (-1)^m$. For $i \in \{m+k+1, m+k+2, ..., n\}$ $\alpha - \lambda_i > 0$, and, consequently, $\prod_{i=m+k+1}^{n}(\alpha - \lambda_i) > 0$. Thus we have proved $\text{sgn}(Q_G(\alpha)) = (-1)^m$. □

Here are some simple and useful consequences of the Lemma 4. In a connected graph $G$, if, $\lambda_2 = \alpha$, then $Q_G(\alpha) < 0$; or, if $\lambda_3 = \alpha$ and $\lambda_2 > \lambda_3$, then $Q_G(\alpha) > 0$.

Now, let us show the main theorem.

**The Generalized RS-theorem.** Let $G$ be the graph in Figure 3, with a cut-vertex $u$. The components of the graph $G$-$u$, the graphs $G_1, ..., G_n$, are connected graphs. For $a > 0$ it holds:

1) If at most one of the graphs $G_1, ..., G_n$ has index $a$, and for the rest of them the indices are less than $a$, then $\lambda_2(G) < a$.
2) If at least two of the graphs $G_1, ..., G_n$ have indices $a$, and for the rest of them the indices are not greater than $a$, then $\lambda_2(G) = a$.
3) If at least two of the graphs $G_1, ..., G_n$ have indices greater than $a$, then $\lambda_2(G) > a$.
4) If only one of the graphs $G_1, ..., G_n$ has index greater than $a$, and at least one of the remaining graphs has index $a$, and for the rest of them the indices are less than $a$, then $\lambda_2(G) > a$.

5) If for one of the graphs $G_1,...,G_n$ $\lambda_1(G_i) > \lambda_2(G_i) > a$ holds, and for the rest of them the indices are less than $a$, then $\lambda_2(G) > a$.

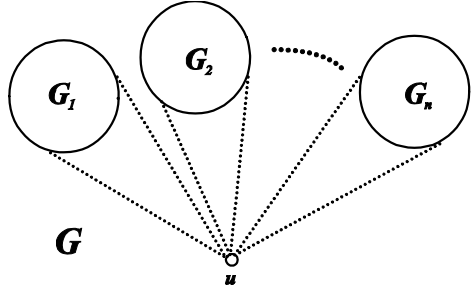

**Figure 3**

**Proof.** 1) If $\lambda_1(G_1), \lambda_1(G_2),...,\lambda_1(G_n) < a$, then, $\lambda_1\left(\bigcup_{i=1}^{n} G_i\right) < a$ and $\lambda_2\left(\bigcup_{i=1}^{n} G_i\right) < a$, and, therefore, by the Interlacing theorem, $\lambda_2(G) < a$.

Now consider the case when exactly one of the graphs $G_i$ has index $a$, say $\lambda_1(G_1) = a$, and $\lambda_1(G_2), \lambda_1(G_3),...,\lambda_1(G_n) < a$. $\lambda_1(G_1) = a$ implies $\lambda_2(G_1) < a$. Then, $\lambda_1\left(\bigcup_{i=1}^{n} G_i\right) = a$ and $\lambda_2\left(\bigcup_{i=1}^{n} G_i\right) < a$, and, therefore, by the Interlacing theorem, $\lambda_2(G) \leq a$. $\lambda_1(G_1) = a$ implies $P_{G_1}(a) = 0$, and $\lambda_1(G_2), \lambda_1(G_3),...,\lambda_1(G_n) < a$ implies $P_{G_2}(a), P_{G_3}(a),...,P_{G_n}(a) > 0$. Applying the Schwenk's lemma at the vertex $u$ of the graph $G$, we compute the characteristic polynomial of $G$.

$$P_G(a) = aP_{G_1}(a)\cdot...\cdot P_{G_n}(a) -$$

$$-\left(\sum_{v \in Adju \cap G_1} P_{G_1-v}(a) + 2\sum_{C \in C(u) \cap G_1} P_{G_1-V(C)}(a)\right) P_{G_2}(a)\cdot...\cdot P_{G_n}(a) -$$

$$-P_{G_1}(a)\left(\sum_{v \in Adju \cap (G-G_1)} P_{G-G_1-u-v}(a) + 2\sum_{C \in C(u) \cap (G-G_1)} P_{G-G_1-V(C)}(a)\right)$$

$$= -\left(\sum_{v \in Adju \cap G_1} P_{G_1-v}(a) + 2\sum_{C \in C(u) \cap G_1} P_{G_1-V(C)}(a)\right) P_{G_2}(a)\cdot...\cdot P_{G_n}(a)$$

Let us introduce the graph $K = G - (G_2 \cup ... \cup G_n)$. Then,

$$P_K(a) = aP_{G_1}(a) - \sum_{v \in Adju \cap G_1} P_{G_1-v}(a) - 2\sum_{C \in C(u) \cap G_1} P_{G_1-V(C)}(a)$$

$$= -\sum_{v \in Adju \cap G_1} P_{G_1-v}(a) - 2\sum_{C \in C(u) \cap G_1} P_{G_1-V(C)}(a).$$

Then, $P_G(a) = P_K(a) \cdot P_{G_2}(a) \cdot \ldots \cdot P_{G_n}(a)$. $P_K(a) < 0$ by the Lemma 1, and therefore, $P_G(a) < 0$, implying $\lambda_2(G) < a$.

2) Let $\lambda_1(G_1) = a$, $\lambda_1(G_2) = a$, and $\lambda_1(G_3), \lambda_1(G_4), \ldots, \lambda_1(G_n) \leq a$. Then, $\lambda_1\left(\bigcup_{i=1}^n G_i\right) = a$ and $\lambda_2\left(\bigcup_{i=1}^n G_i\right) = a$, and, therefore, $\lambda_1(G) > a$ and $\lambda_2(G) = a$, by the Interlacing theorem.

3) If indices of at least two graphs $G_1, \ldots, G_n$ are greater than $a$, say, $\lambda_1(G_1) > a$ and $\lambda_1(G_2) > a$, then $\lambda_1\left(\bigcup_{i=1}^n G_i\right) > a$ and $\lambda_2\left(\bigcup_{i=1}^n G_i\right) > a$, and, therefore, $\lambda_2(G) > a$, by the Interlacing theorem.

4) If the index of one of the graphs $G_1, \ldots, G_n$ is greater than $a$, say $\lambda_1(G_1) > a$, and one is equal to $a$, say $\lambda_1(G_2) = a$, and the remaining indices are less than $a$, we will prove that $\lambda_2(G) > a$. First, let us notice that it is sufficient to prove that $\lambda_2(G) > a$ holds in case when graph $G - u$ consists only of the two mentioned components $G_1$ and $G_2$. The Interlacing theorem implies $\lambda_2(G) \geq a$.

If $\lambda_2(G_1) > a$, then $\lambda_2(G_1 \cup G_2) > a$, and, therefore, $\lambda_2(G) > a$.

Now, let us consider the case $\lambda_2(G_1) = a$.

Let us introduce the graphs $H = G - G_1$ and $K = G - G_2$. Applying the Schwenk lemma at the vertex $u$ of the graphs $H$ and $K$ we get:

$$P_K(\lambda) = \lambda P_{G_1}(\lambda) - \sum_{v \in Adju \cap G_1} P_{G_1 - v}(\lambda) - 2 \sum_{C \in C(u) \cap K} P_{K-C}(\lambda)$$
$$P_H(\lambda) = \lambda P_{G_2}(\lambda) - \sum_{v \in Adju \cap G_2} P_{G_2 - v}(\lambda) - 2 \sum_{C \in C(u) \cap H} P_{H-C}(\lambda)$$

For the graph $G$ we get:

$$P_G(\lambda) = \lambda P_{G_1}(\lambda) P_{G_2}(\lambda) -$$
$$- \left( \sum_{v \in Adju \cap G_1} P_{G_1 - v}(\lambda) + 2 \sum_{C \in C(u) \cap K} P_{K-C}(\lambda) \right) P_{G_2}(\lambda) -$$
$$- \left( \sum_{v \in Adju \cap G_2} P_{G_2 - v}(\lambda) + 2 \sum_{C \in C(u) \cap H} P_{H-C}(\lambda) \right) P_{G_1}(\lambda)$$

Finally, $P_G(\lambda) = P_K(\lambda) P_{G_2}(\lambda) + P_{G_1}(\lambda) P_H(\lambda) - \lambda P_{G_1}(\lambda) P_{G_2}(\lambda)$.

Since $P_{G_2}(a) = 0$, we have $P_G(a) = P_{G_1}(a) P_H(a)$. Since $\lambda_1(G_2) = a$, Lemma 1 implies $P_H(a) < 0$.

a) If $\lambda_2(G_1) < a < \lambda_1(G_1)$, then $P_{G_1}(a) < 0$, implying $P_G(a) > 0$, and, therefore, $\lambda_2(G) > a$.

b) Let us consider the case $\lambda_2(G_1) = \lambda_3(G_1) = ... = \lambda_k(G_1) = a$ and $\lambda_{k+1}(G_1) < a$, $k \geq 2$.

By the Interlacing theorem we have $\lambda_2(K) \geq a$.

If $\lambda_2(K) > a$, then $\lambda_2(G) > a$, also by the Interlacing theorem.

Now, let us consider the case $\lambda_2(K) = a$. Lemma 2 implies $\lambda_3(K) = \lambda_4(K) = ... = \lambda_k(K) = a$. Knowing this, let us introduce the polynomial $Q_K(\lambda)$ by $P_K(\lambda) = (\lambda - a)^{k-1} Q_K(\lambda)$. Similarly, we have $P_{G_1}(\lambda) = (\lambda - a)^{k-1} Q_{G_1}(\lambda)$ and $P_{G_2}(\lambda) = (\lambda - a) Q_{G_2}(\lambda)$. Notice that, by Lemma 4, $Q_{G_1}(a) < 0$. Now,

$$P_G(\lambda) = (\lambda - a)^{k-1} \cdot Q_K(\lambda) \cdot (\lambda - a) \cdot Q_{G_2}(\lambda) +$$
$$+ (\lambda - a)^{k-1} \cdot Q_{G_1}(\lambda) \cdot P_H(\lambda) -$$
$$- \lambda (\lambda - a)^{k-1} \cdot Q_{G_1}(\lambda) \cdot (\lambda - a) \cdot Q_{G_2}(\lambda)$$

Let us introduce the polynomial $Q_G(\lambda)$ by $P_G(\lambda) = (\lambda - a)^{k-1} Q_G(\lambda)$. Then, $Q_G(a) = Q_{G_1}(a) \cdot P_H(a)$. Since, $Q_{G_1}(a) < 0$ and $P_H(a) < 0$, we get $Q_G(a) > 0$, implying $\lambda_2(G) > a$.

5) If $\lambda_1(G_i) > \lambda_2(G_i) > a$ for some $i$, then $\lambda_2(G) > a$, by the Interlacing theorem. □

The following statement is the direct consequence of the Generalized *RS*-theorem.

**Corollary 1.** Let $G$ be a connected graph with the cut-point $u$, as in Figure 3. Let $\alpha_1$ and $\alpha_2$, $\alpha_1 \geq \alpha_2$, be the two largest indices of the connected components of the graph *G-u*. Then:

1) If $\alpha_1 = \alpha_2$, then $\lambda_2(G) = \alpha_1$.

2) If $\alpha_1 > \alpha_2$, then $\alpha_2 < \lambda_2(G) < \alpha_1$.

**Proof.** We prove the first statement by putting $\alpha_1 = \alpha_2 = a$ and applying the *GRS*-theorem. To prove the second statement, we put $\alpha_2 = a$, and apply the theorem. □

## 4. Remarks and examples of application

The Generalized *RS*-theorem does not provide the information on $\lambda_2(G)$ if $\lambda_1(G_i) > a$ and $\lambda_2(G_i) \leq a$, for one of the graphs $G_1,...,G_n$, and for the rest of them the indices are less than $a$. In that case it is possible that $\lambda_2(G) < a$, $\lambda_2(G) > a$ or $\lambda_2(G) = a$. We will show this in the Examples 1 and 2.

**Example 1.** Let us consider graphs $M$ and $N$ shown in Figure 4. In both graphs, after removal of the cut-vertex $u$, we get the same two connected graphs $G_1$ and $G_2$ for which $\lambda_1(G_1) > \sqrt{3}$ and $\lambda_1(G_2) < \sqrt{3}$. For the second largest eigenvalue of the graphs $M$ and $N$, we get $\lambda_2(M) < \sqrt{3}$, and $\lambda_2(N) > \sqrt{3}$.

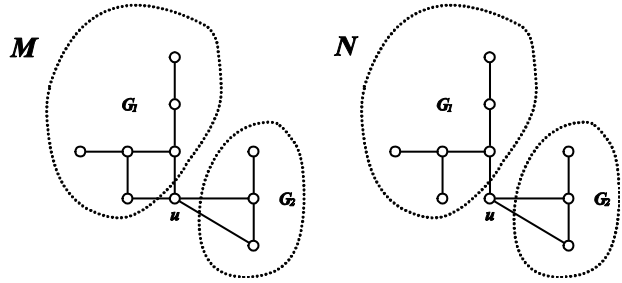

Figure 4

**Example 2.** Let us consider the graph $G$ in Figure 5. After removal of the cut-vertex $u$, we get three connected graphs: $G_1$, $G_2$ and $G_3$ for whose indices it holds: $\lambda_1(G_1) > \sqrt{3}$, $\lambda_1(G_2) < \sqrt{3}$ and $\lambda_1(G_3) < \sqrt{3}$. For the original graph $G$ the following is true: $\lambda_2(G) = \lambda_3(G) = \sqrt{3}$.

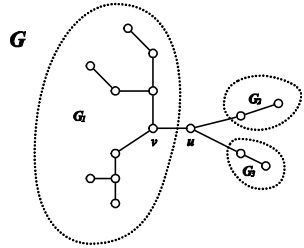

Figure 5

Note that for the construction of these graphs we used the fact that we know the only two trees with the property $\lambda_1 = \sqrt{3}$. These trees can be seen in graph $G_1$ (Figure 5) after removal of the vertex $v$. Note also, that the graph $G_1$ itself is constructed using the Generalized *RS*-theorem, i.e. if these two trees, with the property $\lambda_1 = \sqrt{3}$, are connected via the vertex $v$ as in Figure 5, then the Generalized *RS*-theorem states that $\lambda_2(G_1) = \sqrt{3}$ (and, also, $\lambda_1(G_1) > \sqrt{3}$).

Also, looking again at the graph $G$, we see that, after removal of the vertex $v$, we get three graphs with the index $\sqrt{3}$.

Now, let us show an application of the Generalized *RS*-theorem. Besides determining whether the second largest eigenvalue of a graph is greater than, smaller than, or equal to some value, we can use this theorem to find an upper and a lower bound for the second largest eigenvalue of a graph.

**Example 3.** Consider the graph $G$ in Figure 6. After removal of the cut-vertex $u$, we get three components. Two of them are supergraphs of the graphs with the index $\sqrt{3}$ (and the third one is a subgraph), and, therefore $\lambda_2(G) > \sqrt{3}$, by the Generalized *RS*-theorem. These three components are also subgraphs of Smith graphs (their indices are smaller than 2), and so we conclude that $\lambda_2(G) < 2$. Thus we have found the bounds for the second largest eigenvalue of the graph $G$: $\sqrt{3} < \lambda_2(G) < 2$.

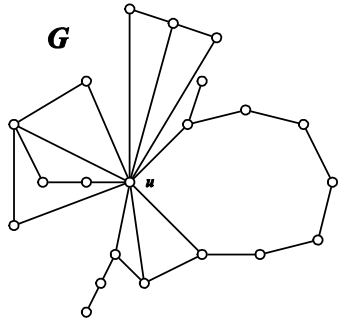

Figure 6

The tables of spectra and other relevant data for connected graphs with the small number of vertices are available in the literature. For example, in [4], there are spectra of all connected graphs with up to six vertices, all trees with up to nine vertices and all cubic graphs with up to twelve vertices. Therefore, in many cases after removing the cut vertex $u$ of the graph $G$, we know the indices of the connected components of $G$-$u$ without computing, and consequently, the bounds for the second largest eigenvalue.

**Example 4.** Consider the graph in Figure 7. Graph $G_1$ is a cubic graph and its index is 3. Index of the graph $G_3$ is $2\sqrt{2}$ ($\approx 2.82843$). Graph $G_2$ is a path and its index is smaller than 2. The indices of the graphs $G_4$ and $G_5$ are 2.1010 and 2.0421 respectively. Using the Corollary 1, we conclude that $2\sqrt{2} < \lambda_2(G) < 3$. If we actually compute $\lambda_2$, we get: $\lambda_2(G) \approx 2.9365$.

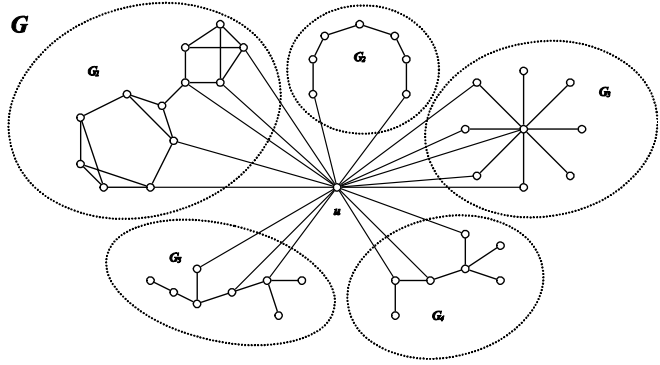

Figure 7

## 5. Conclusion

In our previous research the *RS*-theorem was a valuable tool in search for the reflexive graphs. Given graphs with the index *2* (Smith graphs), and using the *RS*-theorem, we were able to establish whether the second largest eigenvalue of a given graph is less than, greater than or equal to *2*, in many cases of connected graphs with a cut-point. This is done in the following way. We remove the cut-vertex *u* of the graph *G* and, then, compare the components of the graph *G-u* to Smith graphs. We examine whether they are Smith graphs, their subgraphs or supergraphs, and see if the *RS*-theorem could be applied.

In this paper we have stated and proved the Generalized *RS*-theorem. Within some classes of connected graphs with a cut-vertex the Generalized *RS*-theorem gives us the answer whether $\lambda_2(G) < a$, $\lambda_2(G) > a$ or $\lambda_2(G) = a$, for some $a > 0$. The method is similar to the one described in the previous paragraph. Given the graph G with a cut-vertex *u,* we use the indices of the components of the graph *G-u,* compare them to *a*, and make conclusions about the second largest eigenvalue of *G*.

Also, we believe that Lemmas 1-4, used in the proof of the Generalized *RS*-theorem, will be useful tools in our future research. For example, many times we have had the situation that after computing the characteristic polynomial, we get $P_G(a) = 0$, and thus we know that *a* is an eigenvalue of *G*, but we do not know whether it is, say, $\lambda_2$ or $\lambda_3$. Lemma 4 is an effective algebraic tool in such cases. Finally, using the Corollary 1 we can determine the bounds for the second largest eigenvalue for many graphs with a cut-vertex.

**Acknowledgement**. The programming package *newGRAPH* [23] was used for computations of eigenvalues of graphs.

M. Rašajski is grateful to the Serbian Ministry of Science and Technological Development for financial support through the Grant No. 174033.